\documentclass{autart}
\usepackage{amsmath,amssymb,amsfonts}
\usepackage{algorithmic}
\usepackage{graphicx}
\usepackage{textcomp}
\usepackage{mathrsfs}
\usepackage{bm}
\usepackage{xcolor}
\usepackage{ntheorem}
\newtheorem*{Proof}{Proof}       
\newtheorem{Proposition}{\textbf{Proposition}} 
\newtheorem{Notation}{\textbf{Notation}} 
\newtheorem{Definition}{\textbf{Definition}}
\newtheorem{Lemma}{\textbf{Lemma}}
\newtheorem{Theorem}{\textbf{Theorem}} 
\newtheorem{Remark}{\textbf{Remark}} 
\allowdisplaybreaks[4]
\usepackage{epsfig}
\usepackage{subfig}
\usepackage[title]{appendix}
\usepackage{float} 
\graphicspath{{pictures/}}

\allowdisplaybreaks

\begin{document}

\begin{frontmatter}

\title{Bilateral boundary control of an input delayed \\ 2-D reaction-diffusion equation \thanksref{t1}\thanksref{t2}} 
\thanks[t1]{This paper was not presented at any conference.}
\thanks[t2]{This work is supported by National Natural Science Foundation of China (62173084, 61773112) and the Project of Science and Technology Commission of Shanghai Municipality, China (23ZR1401800, 22JC1401403), in part by the Fundamental Research Funds for the Central Universities and Graduate Student Innovation Fund of Donghua University CUSF-DH-D-2023043.}

\author[engineer]{Dandan Guan}\ead{gdd\_dhu@163.com},    
\author[engineer]{Yanmei Chen}\ead{xiaobanshiguang@163.com},               
\author[engineer,engineer2]{Jie Qi\thanksref{cor}}\ead{jieqi@dhu.edu.cn}\thanks[cor]{Corresponding author}, 
\author[mathematics]{Linglong Du}\ead{matdl@dhu.edu.cn}  
\address[engineer]{College of Information Science and Technology, Donghua University, Shanghai, 201620, China}                                                
\address[engineer2]{Engineering Research Center of Digitized Textile and Apparel Technology, Ministry of Education, Donghua University
	Shanghai, 201620, China}            
\address[mathematics]{Department of Applied Mathematics, Donghua University, Shanghai, 201620, China}        

\begin{keyword}                           
	2-$D$ Reaction-Diffusion; Input delay;  Backstepping;  Bilateral boundary control.             
\end{keyword}                             

\begin{abstract}                          
	In this paper, a delay compensation design method based on PDE backstepping is developed for a two-dimensional reaction-diffusion partial differential equation (PDE) with bilateral input delays. The PDE is defined in a rectangular domain, and the bilateral control is imposed on a pair of opposite sides of the rectangle. To represent the delayed bilateral inputs, we introduce two 2-D transport PDEs that form a cascade system with the original PDE.  A novel set of backstepping transformations is proposed for delay compensator design, including one Volterra integral transformation and two affine Volterra integral transformations.  Unlike the kernel equation for 1-D PDE systems with delayed boundary input, the resulting kernel equations for the 2-D system have singular initial conditions governed by the Dirac Delta function.  Consequently, the kernel solutions are written as a double trigonometric series with singularities.   To address the challenge of stability analysis posed by the singularities, we prove a set of inequalities by using the Cauchy-Schwarz inequality, the 2-D Fourier series, and the Parseval's theorem.  A numerical simulation illustrates the effectiveness of the proposed delay-compensation method.
	
\end{abstract}

\end{frontmatter}

\section{Introduction}
Reaction-diffusion equations are usually used to describe the density concentration of a substance, a population change in time and space, governed by Fick's law and reaction factors, such as chemical reaction, birth, and death. Many of these processes are unstable due to the reaction terms, naturally applied in chemical reactions \cite{orlov2002discontinuous}, thermal fluids \cite{eleiwi2018observer}, population dynamics \cite{kondo2010reaction}, etc., which gives rise to various challenges for control design, especially for the problem with delay which would greatly increase the difficulty of stabilizing the unstable system \cite{zhang2006optimal}.

In past decades, the delay-compensated control for the reaction-diffusion equation has been well developed. An early work is \cite{krstic2009control}, which considers the stabilization of a reaction-diffusion PDE with arbitrary long input delay by the PDE backstepping boundary controller.  \cite{steeves2021input,wang2021delay,qi2019stabilization,qi2020compensation,hashimoto2016stabilization,wang2017output}  extend the method of \cite{krstic2009control} to solve different delay-compensated control problems for  reaction-diffusion systems, in which,  \cite{steeves2021input} studies the problem of prescribed–time stabilization for a reaction-diffusion PDE with constant input delay.  In \cite{qi2019stabilization}, a delay compensator is designed for distributed in-domain control.   A state delay problem is considered  in \cite{hashimoto2016stabilization}, 
and solved by a backstepping controller.  The spatially-varying delay compensator is designed in \cite{qi2020compensation}, while an unknown input delayed control problem is solved by a delay-adaptive predictor controller in  \cite{wang2021delay}. The output regulator delay is worked out in \cite{wang2017output}.
By a similar approach, \cite{koga2020delay} presents the control design for the Stefan problem under actuator delay. Different from the above papers,  \cite{fridman2009exponential} applies Lyapunov-Krasovskii functional methods and a linear operator inequality to solve the delay compensation problem and achieve the exponential stability of the system with time-varying delays. Another paper using Lyapunov-Krasovskii functional methods is \cite{9732660}, in which the authors design a new sub-predictor-based T-S fuzzy control law that can achieve compensation for large delay via using the proposed observer.  Besides, \cite{prieur2018feedback} designs a stabilizing feedback boundary control for a reaction-diffusion PDE with input delay by using the spectral decomposition method and Artstein transformation, whereas \cite{lhachemi2020pi} applies a similar method to design PI regulation of a reaction-diffusion PDE with delay. 
Another way to compensate for the input delay is predictor feedback control which stabilizes systems via predicting the future value of the state \cite{karafyllis2011stabilization,selivanov2018boundary}. \cite{karafyllis2011stabilization} proposes a global stabilization strategy for the nonlinear system with delayed input, 
and
\cite{selivanov2018boundary} handles the reaction-diffusion PDE with delay and data sampling based on observer. A similar study is reported in \cite{katz2021delayed}, the Lyapunov functions and Halanays inequality are applied for finite-dimensional observer-based control of a reaction-diffusion PDE with fast-varying unknown input and known output delays.

The bilateral control can reduce control efforts compared to the unilateral control, especially for the system with large unstable coefficients \cite{vazquez2016bilateral}, so more attention has been paid to bilateral control, such as \cite{liu2020boundary,chen2021folding,bekiaris2019nonlinear}. 
In \cite{bekiaris2019nonlinear}, the authors employ the PDE backstepping method to design a nonlinear bilateral full-state feedback controller for a class of viscous Hamilton-Jacobi PDEs. Another important result is reported in \cite{liu2020boundary}, bilateral backstepping boundary controller is developed for a reaction-diffusion PDE  in $n$-D  domain symmetric about the $x$-axis. Based on state feedback and state estimation, \cite{chen2021folding} uses the backstepping method and fold approach to construct the bilateral output-feedback controller for an unstable 1-D reaction-diffusion system with spatially-varying reaction. In addition, bilateral control has also been applied into practice, for example,  traffic flow control \cite{yu2019bilateral}, 
the tracking adaptive control issue of the master-slave flexible system \cite{wang2021bilateral}.

Although high-dimensional models have more applications, such as current reaction in tokamak device \cite{mavkov2018experimental} and temperature control of  screw extruders \cite{lipar2013extruder},  fewer research work addresses delayed PDEs in the high-dimensional domain. An example is  \cite{selivanov2019delayed}, which proposes a $H_\infty$ robust controller for 2-D diffusion systems with small delayed pointlike measurements and small delayed input by using the Lyapunov-Krasovskii functional approach.  For a particular cylindrical surface,  \cite{qi2019control} stabilizes an unstable  2-D reaction-diffusion equation with input delay by the PDE backstepping. In the paper, we consider a bilateral input delay problem of a reaction-diffusion system in a rectangular domain. An example application of the considered model is extruder temperature control by heating/cooling cylindrical barrels. Since the temperature is symmetrical with the central axis of the cylinder, the temperature distribution is the same for each section crossing the central axis. Hence, the PDE domain defined in the cylinder can be simplified to a 2-D rectangle. Owing to two actuators being subject to delays, the classical backstepping transformation \cite{krstic2009control} is no longer applicable. We propose a novel set of backstepping transformations containing a Volterra  transformation  \cite{vazquez2016bilateral,liu2020boundary} and two affine Volterra  transformations. Since the 4-D kernel functions in the affine Volterra transformations depend on the 2-D kernel function in the Volterra transformation at one of their boundaries, a Dirac Delta function is employed to build this relationship by reducing one dimension.  
Thus,  the resulting kernel equations involve Dirac delta function as their initial conditions. We solve the kernel equations, express the solution as double series, and prove the series converges, including the convergence of the series at singularity $s=0, ~x=\theta$ to the initial conditions governed by Dirac delta function. The singularity in 2-D domain poses more challenges for the proof of the operator boundedness than that in 1-D because it needs to compute the $H^1$ norm of the delayed actuator state in 2-D domain, rather than the $L^2$ norm in 1-D domain. We solve the difficulty by proving a group of inequalities containing integration of kernel functions by using 2-D Fourier series \cite{haberman2012applied}, Cauchy-Schwarz inequality and  Parseval's theorem, and then obtain the exponential stability of the closed-loop system in terms of $L^2$ norm for the state and $H^1$ norm for the actuator states. Finally, the theoretical results are demonstrated on a 2-D square domain by simulation. 


The contribution of this paper is summarized as:
\begin{itemize}
\item First, we propose a novel set of backstepping transformations to solve the bilateral control problem of a 2-D reaction-diffusion system in presence of input delays for both actuators. 
\item To address the  stability analysis challenge posed by the singularity, we prove a group of inequalities containing integration of kernel functions by using Cauchy-Schwarz inequality,  2-D Fourier series, and  Parseval's theorem.
\end{itemize}

This paper is organized as follows. Section 2 presents the control design for the system with input delay. The stability analysis is shown in Section 3 and supportive simulation results are provided in Section 4. The paper ends with concluding remarks and future works presented in Section 5.
\setlength{\abovedisplayskip}{3pt}
\setlength{\belowdisplayskip}{3pt}
\begin{Notation}\rm
In this paper, the following notations are used to define the domains: 
\begin{align*}
	\mathcal{D}_1&=\{(x,y)|-L\le x\le L, -l\le y\le l\},\\
	\mathcal{D}_2&=\{(x,s)|-L\le x\le L, 0\le s\le \tau\},\\
	\Gamma_1&=\{(y,\xi)|-l\le y\le l, -y\le \xi\le y\},\\
	\Gamma_2&=\{(x,\theta,s,\xi)|-L\le x, \theta\le L, 0\le s\le \tau, -y\le\xi\le y\}.
\end{align*} 
We also define the norms for $f(x,y)\in L^2(\mathcal{D}_1)$,
\begin{align*}
	\Vert f\Vert_{L^2}^2&=\int_{-L}^{L}\int_{-l}^{l}f^2(x,y)dydx,\\
	\Vert f\Vert_{H^1}^2&=\Vert f\Vert_{L^2}^2+\Vert f_x\Vert_{L^2}^2+\Vert f_y\Vert_{L^2}^2,
\end{align*}
and for $g(x,s)\in L^2(\mathcal{D}_2)$,
\begin{align*}
	\Vert g\Vert_{L^2}^2&=\int_{-L}^{L}\int_{0}^{\tau}g^2(x,s)dsdx,\\
	\Vert g\Vert_{H^1}^2&=\Vert g\Vert_{L^2}^2+\Vert g_x\Vert_{L^2}^2+\Vert g_s\Vert_{L^2}^2.
\end{align*}
\end{Notation}
\section{Controller design} 
\subsection{Problem description}
Consider the following unstable $2$-$D$ reaction-diffusion system with input delay $\tau>0$ on opposite boundaries $y=l$ and $y=-l$,
\begin{align}
u_t(x,y,t)&=u_{xx}(x,y,t)+u_{yy}(x,y,t)+\lambda u(x,y,t)\label{original-u},\\
u(L,y,t)&=u(-L,y,t)=0\label{boundary-u-o},\\
u(x,l,t)&=U_1(x,t-\tau)\label{control-U-1},\\
u(x,-l,t)&=U_2(x,t-\tau)\label{control-U-2},\\
u(x,y,0)&=u_0(x,y),\label{initial-u}
\end{align}
where  $(x,y)\in\mathcal{D}_1$, $t>0$. The reaction coefficient $\lambda>0$. $U_1(x,t-\tau)$ and $U_2(x,t-\tau)$ are the control inputs that will be determined later. 


Introduce two 2-D transport PDEs representation of the delay on the distributed inputs, which results in the following equivalent cascade system:
\begin{align}
u_t(x,y,t)&=u_{xx}(x,y,t)+u_{yy}(x,y,t)+\lambda u(x,y,t)\label{cascade-u},\\
u(-L,y,t)&=u(L,y,t)=0,\label{cascade-boundary-u}\\
u(x,l,t)&=v_1(x,0,t)\label{cascade-v1},\\
u(x,-l,t)&=v_2(x,0,t)\label{cascade-v2},\\
u(x,y,0)&=u_0(x,y)\label{cascade-initial-u0},\\
\partial_t v_1(x,s,t)&=\partial_s v_1(x,s,t)\label{cascade-v1-pde}, \\
\partial_t v_2(x,s,t)&=\partial_s v_2(x,s,t)\label{cascade-v2-pde},\\
v_1(x,\tau,t)&=U_1(x,t),\label{cascade-v-control1}\\
v_2(x,\tau,t)&=U_2(x,t)\label{cascade-v-control2},\\
v_1(x,s,0)&=v_{10}(x,s),\label{cascade-v-initial1}\\
v_2(x,s,0)&=v_{20}(x,s)\label{cascade-v-initial2},
\end{align}
where $(x,s)\in\mathcal{D}_2$ for states $v_1(\cdot,\cdot,t)$ and  $v_2(\cdot,\cdot,t)$.

\subsection{The backstepping transformation}
Apply the PDE backstepping method for control
design. First, introduce a stable target system:
\begin{align}
w_t(x,y,t)&=w_{xx}(x,y,t)+w_{yy}(x,y,t),\label{target-w}\\
w(-L,y,t)&=w(L,y,t)=0\label{target-w-boundary},\\
w(x,l,t)&=z_1(x,0,t)\label{target-z1-control},\\
w(x,-l,t)&=z_2(x,0,t)\label{target-z2-control},\\
w(x,y,0)&=w_0(x,y)\label{terget-w-initial},\\
\partial_t z_1(x,s,t)&=\partial_s z_1(x,s,t)\label{target-z1-pde},\\
\partial_t z_2(x,s,t)&=\partial_s z_2(x,s,t)\label{target-z2-pde},\\
z_1(x,\tau,t)&=0\label{target-z-boundary1},\\
z_2(x,\tau,t)&=0\label{target-z-boundary2},\\
z_1(x,s,0)&=z_{10}(x,s),\label{target-z-initial1}\\
z_2(x,s,0)&=z_{20}(x,s)\label{target-z-initial2},
\end{align}
where the transport systems $z_1$ and $z_2$ have mild solutions:
\begin{align}\label{mild-solution-zi}
z_i(x,s,t)=\left\{
\begin{aligned}
	&z_{i0}(x,s+t),&0<s+t\le \tau\\
	&0,&\tau < s+t
\end{aligned}\
,~~i=1,2\right.,
\end{align}
which implies $z_i(x,s,t)$ is stable in finite time.    To transform the original cascade system \eqref{cascade-u}-\eqref{cascade-v-initial2} into  \eqref{target-w}-\eqref{target-z-initial2}, we propose a set of backstepping integral transformations as follows:
\begin{align}
&w(x,y,t)=u(x,y,t)-\int_{-y}^{y}p(y,\xi)u(x,\xi,t)d\xi,\label{trans-w}\\
&z_1(x,s,t)=v_1(x,s,t)\nonumber\\
&-\int_{-L}^{L}\int_{-l}^{l}\gamma_1(x,\theta,s,\xi)u(\theta,\xi,t)d\xi d\theta\nonumber\\
&-\int_{-L}^{L}\int_{0}^{s}\partial_\xi\gamma_1(x,\theta,s-r,l)v_1(\theta,r,t)dr d\theta\nonumber\\
&+\int_{-L}^{L}\int_{0}^{s}\partial_\xi\gamma_1(x,\theta,s-r,-l)v_2(\theta,r,t)dr d\theta\label{trans-z1},\\
&z_2(x,s,t)=v_2(x,s,t)\nonumber\\
&-\int_{-L}^{L}\int_{-l}^{l}\gamma_2(x,\theta,s,\xi)u(\theta,\xi,t)d\xi d\theta\nonumber\\
&-\int_{-L}^{L}\int_{0}^{s}\partial_\xi\gamma_2(x,\theta,s-r,l)v_1(\theta,r,t)dr d\theta\nonumber\\
&+\int_{-L}^{L}\int_{0}^{s}\partial_\xi\gamma_2(x,\theta,s-r,-l)v_2(\theta,r,t)dr d\theta\label{trans-z2},
\end{align}
where $\partial_\xi\gamma_i (\cdot,\cdot,\cdot,l)=\frac{\partial \gamma_i(\cdot,\cdot,\cdot,\xi)}{\partial \xi}|_{\xi=l}$, $i=1,2$,  transformation \eqref{trans-w} is the Volterra transformation,  transformations \eqref{trans-z1} and \eqref{trans-z2} are the affine Volterra transformations. The kernel function $p(y,\xi)$ is defined in domain $\Gamma_1$,  and  $\gamma_1$ and $\gamma_2$ are both defined in $\Gamma_2$.

According to \cite{vazquez2016bilateral,liu2020boundary}, the kernel function $p(y,\xi)$  satisfies the following PDE:
\begin{align}
p_{yy}(y,\xi)-p_{\xi\xi}(y,\xi)&=\lambda p(y,\xi),\label{main-p}\\
p(y,y)&=-\frac{\lambda}{2}y,\label{p-bnd1}\\
p(y,-y)&=0,\label{p-bnd2}
\end{align}
whose solution is written as:
\begin{align}\label{solut-p}
p(y,\xi)=-\frac{sgn(y)}{2}\sqrt{\lambda\frac{y+\xi}{y-\xi}}I_1\left(\sqrt{\lambda(y^2-\xi^2)}\right),
\end{align}
where $I_1(\cdot)$ is the first-order modified Bessel function of the first kind. The approach to solving the equation \eqref{main-p}-\eqref{p-bnd2} can be found in \cite{vazquez2016bilateral,vazquez2016explicit} .

After a lengthy computation based on the equivalent relation between system \eqref{cascade-u}-\eqref{cascade-v-initial2} and system \eqref{target-w}-\eqref{target-z-initial2}, one can get the equations of kernels $\gamma_i$, $i=1,2$:
\begin{align}
\partial_s\gamma_1(x,\theta,s,\xi)=&\partial_{\xi\xi}\gamma_1(x,\theta,s,\xi)+\partial_{\theta\theta}\gamma_1(x,\theta,s,\xi)\nonumber\\
&+\lambda\gamma_1(x,\theta,s,\xi),\label{equ-gamma1}\\
\gamma_1(x,L,s,\xi)=&\gamma_1(x,-L,s,\xi)=0,\label{boun1-gamma1}\\
\gamma_1(x,\theta,s,l)=&\gamma_1(x,\theta,s,-l)=0,\label{boun2-gamma1}\\
\gamma_1(x,\theta,0,\xi)=&p(l,\xi)\delta(x-\theta),\label{initial-gamma1}\\
\partial_s\gamma_2(x,\theta,s,\xi)=&\partial_{\xi\xi}\gamma_2(x,\theta,s,\xi)+\partial_{\theta\theta}\gamma_2(x,\theta,s,\xi)\nonumber\\
&+\lambda\gamma_2(x,\theta,s,\xi),\label{equ-gamma2}\\
\gamma_2(x,L,s,\xi)=&\gamma_2(x,-L,s,\xi)=0,\label{boun1-gamma2}\\
\gamma_2(x,\theta,s,l)=&\gamma_2(x,\theta,s,-l)=0,\label{boun2-gamma2}\\
\gamma_2(x,\theta,0,\xi)=&-p(-l,\xi)\delta(x-\theta),\label{initial-gamma2}
\end{align}
where $\delta(\cdot)$ is the Dirac delta function.

Applying the method of separation of variables, one gets the solutions  of kernel functions $\gamma_1$ and $\gamma_2$:
\begin{align}
\gamma_1(x,\theta,s,\xi)=&\frac{1}{L}\sum_{n=1}^{\infty}\sum_{m=1}^{\infty}e^{\left(-\frac{(n\pi)^2}{4L^2}-\frac{(m\pi)^2}{4l^2}+\lambda\right)s}\nonumber\\
&\cdot p_{1m}\varphi_n(x)\varphi_n(\theta)\phi_m(\xi),\label{solution-gamma1}\\
\gamma_2(x,\theta,s,\xi)=&-\frac{1}{L}\sum_{n=1}^{\infty}\sum_{m=1}^{\infty}e^{\left(-\frac{(n\pi)^2}{4L^2}-\frac{(m\pi)^2}{4l^2}+\lambda\right)s}\nonumber\\
&\cdot p_{2m}\varphi_n(x)\varphi_n(\theta)\phi_m(\xi),\label{solution-gamma2}
\end{align}
where
\begin{align}
&\mathrm{for}~n=1,2,\cdots,~~\varphi_n(\cdot)=\sin\left(\frac{n\pi}{2L}(\cdot+L)\right),\label{varphi}\\
&\mathrm{for}~m=1,2,\cdots,~\phi_m(\xi)=\sin\left(\frac{m\pi}{2l}(\xi+l)\right),\label{phi}\\
&~~~p_{im}=\frac{1}{l}\int_{-l}^{l}\phi_m(\xi)p((-1)^{i-1}l,\xi)d\xi,~i=1,2 \label{p-fourier}.
\end{align} 


\begin{Remark}\rm
It is easy to prove that  kernels  $\gamma_1$ and $\gamma_2$ are convergent as $s>0$. 
As $s=0$, take $\gamma_1$ as an example,
\begin{align*}
	\gamma_1(x,\theta,0,\xi)=&\frac{1}{L}\sum_{n=1}^{\infty}\sum_{m=1}^{\infty} p_{1m}\varphi_n(x)\varphi_n(\theta)\phi_m(\xi)\\
	=&p(l,\xi)\sum_{n=1}^{\infty}   \frac{1}{L} \varphi_n(x)\varphi_n(\theta).
\end{align*}
It is well known \cite{tang2007mathematical} that the Dirac Delta function can be expressed by 
$ \delta(x-\theta)=\sum_{n=1}^{\infty}   \frac{1}{L} \varphi_n(x)\varphi_n(\theta)$, which gives the initial condition \eqref{initial-gamma1}. Consequently, kernels $\gamma_1$ and $\gamma_2$ have singularities at $s=0$ and $x=\theta$ which are  governed by the Delta function. 
\end{Remark}

Combining the boundary conditions \eqref{cascade-v-control1}, \eqref{cascade-v-control2} and \eqref{target-z-boundary1}, \eqref{target-z-boundary2}, substituting $s=\tau$ into transformations \eqref{trans-z1} and \eqref{trans-z2}, we can get the delay-compensated actuators as follows:
\begin{align}
&U_1(x,t)=\int_{-L}^{L}\int_{-l}^{l}\gamma_1(x,\theta,\tau,\xi)u(\theta,\xi,t)d\xi d\theta\nonumber\\
&+\int_{-L}^{L}\int_{t-\tau}^{t}\partial_\xi\gamma_1(x,\theta,t-r,l)U_1(\theta,r)dr d\theta\nonumber\\
&-\int_{-L}^{L}\int_{t-\tau}^{t}\partial_\xi\gamma_1(x,\theta,t-r,-l)U_2(\theta,r)dr d\theta,\label{controller-U1}\\
&U_2(x,t)=\int_{-L}^{L}\int_{-l}^{l}\gamma_2(x,\theta,\tau,\xi)u(\theta,\xi,t)d\xi d\theta\nonumber\\
&+\int_{-L}^{L}\int_{t-\tau}^{t}\partial_\xi\gamma_2(x,\theta,t-r,l)U_1(\theta,r)dr d\theta\nonumber\\
&-\int_{-L}^{L}\int_{t-\tau}^{t}\partial_\xi\gamma_2(x,\theta,t-r,-l)U_2(\theta,r)dr d\theta\label{controller-U2}.
\end{align}
It shows that control $U_1$ and $U_2$ are coupled by both historical inputs as feedback in each controller.
\subsection{Inverse transformation}
A direct way to show the invertibility of the transformation \eqref{trans-w}-\eqref{trans-z2} is to derive its inverse. First, define the inverse transformation as follows,
\begin{align}
&u(x,y,t)=w(x,y,t)+\int_{-y}^{y}q(y,\xi)w(x,\xi,t)d\xi,\label{trans-u}\\
&v_1(x,s,t)=z_1(x,s,t)\nonumber\\
&+\int_{-L}^{L}\int_{-l}^{l}\eta_1(x,\theta,s,\xi)w(\theta,\xi,t)d\xi d\theta\nonumber\\
&+\int_{-L}^{L}\int_{0}^{s}\partial_\xi\eta_1(x,\theta,s-r,l)z_1(\theta,r,t)dr d\theta\nonumber\\
&-\int_{-L}^{L}\int_{0}^{s}\partial_\xi\eta_1(x,\theta,s-r,-l)z_2(\theta,r,t)dr d\theta\label{trans-v1},\\
&v_2(x,s,t)=z_2(x,s,t)\nonumber\\
&+\int_{-L}^{L}\int_{-l}^{l}\eta_2(x,\theta,s,\xi)w(\theta,\xi,t)d\xi d\theta\nonumber\\
&+\int_{-L}^{L}\int_{0}^{s}\partial_\xi\eta_2(x,\theta,s-r,l)z_1(\theta,r,t)dr d\theta\nonumber\\
&-\int_{-L}^{L}\int_{0}^{s}\partial_\xi\eta_2(x,\theta,s-r,-l)z_2(\theta,r,t)dr d\theta\label{trans-v2},
\end{align}
the kernels $q$, $\eta_1$, and $\eta_2$ are expressed as
\begin{align}
q(y,\xi)=&-\frac{sgn(y)}{2}\sqrt{\lambda\frac{y+\xi}{y-\xi}}J_1\left(\sqrt{\lambda(y^2-\xi^2)}\right),\\
\eta_1(x,\theta,s,\xi)=&\frac{1}{L}\sum_{n=1}^{\infty}\sum_{m=1}^{\infty}e^{\left(-\frac{(n\pi)^2}{4L^2}-\frac{(m\pi)^2}{4l^2}\right)s}\nonumber\\
&\cdot q_{1m}\varphi_n(x)\varphi_n(\theta)\phi_m(\xi),\label{solution-eta1}\\
\eta_2(x,\theta,s,\xi)=&-\frac{1}{L}\sum_{n=1}^{\infty}\sum_{m=1}^{\infty}e^{\left(-\frac{(n\pi)^2}{4L^2}-\frac{(m\pi)^2}{4l^2}\right)s}\nonumber\\
&\cdot q_{2m}\varphi_n(x)\varphi_n(\theta)\phi_m(\xi),\label{solution-eta2}
\end{align}
where $J_1(\cdot)$ is the first-order Bessel function, and $
q_{im}=\frac{1}{l}\int_{-l}^{l}\sin\left(\frac{m\pi}{2l}(\xi+l)\right)q((-1)^{i-1}l,\xi)d\xi$, $i=1,2$.

It is worth discussing the difference in control design  between identical delay values and  distinct delay values. In the following remark, we give some hints for control design when the input delays are different.
\begin{Remark}\rm
Let $\tau_1$ and $\tau_2$ denote  distinct delay values for boundary control $U_1$ and $U_2$, respectively. Without loss of generality, assume that $\tau_1 < \tau_2$. The integral transformations \eqref{trans-w} and \eqref{trans-z1} can be kept unchanged, but the integral interval and the integrated state of \eqref{trans-z2} are modified as follows:
\begin{align}
	&z_2(x,s,t)=v_2(x,s,t)\nonumber\\
	&-\int_{-L}^{L}\int_{-l}^{l}\gamma_2(x,\theta,s,\xi)u(\theta,\xi,t)d\xi d\theta\nonumber\\
	&-\int_{-L}^{L}\int_{0}^{\phi(s)}\partial_\xi\gamma_2(x,\theta,s-r,l)z_1(\theta,r,t)dr d\theta\nonumber\\
	&+\int_{-L}^{L}\int_{0}^{s}\partial_\xi\gamma_2(x,\theta,s-r,-l)v_2(\theta,r,t)dr d\theta\label{trans-z21},
\end{align} 
with \begin{align}\label{phi1}
	\phi(s)=\left\{
	\begin{aligned}
		&s,&s\le \tau_1\\
		&\tau_1,&\tau_1 < s
	\end{aligned}\right..
\end{align}			
Based on the transformation, one can obtain a set of kernel equations in a similar method. Due to the space constraints, we omit the details of the control design for the different input  delays.  
\end{Remark}
\section{Stability analysis}
\begin{Theorem}\label{stability-origianl}
\rm For any initial conditions $(u_0,v_{10},v_{20})\in L^2(\mathcal{D}_1)\times H^1(\mathcal{D}_2)\times H^1(\mathcal{D}_2)$ and the compatible conditions holding
\begin{align}
	&u_0(L,y)=0,~~ u_0(-L,y)=0,\\
	&u_0(x,l)=v_{10}(x,0),~~ u_0(x,-l)=v_{20}(x,0),\\
	&v_{10}(x,\tau)=U_1(x,0),~~ v_{20}(x,\tau)=U_2(x,0),
\end{align}
the  closed-loop system \eqref{cascade-u}-\eqref{cascade-v-initial2} with the delay-compensated actuators 	\eqref{controller-U1} and \eqref{controller-U2} admits a unique solution that satisfies 
\begin{align}
	&\Vert u(t)\Vert^2_{L^2}+\Vert v_1(t)\Vert^2_{H^1}+\Vert v_2(t)\Vert^2_{H^1}\nonumber\\
	\le&\alpha_1e^{-\beta_1t}(\Vert u_0\Vert^2_{L^2}+\Vert v_{10}\Vert^2_{H^1}+\Vert v_{20}\Vert^2_{H^1}),
\end{align}
for constants $\alpha_1, \beta_1>0$.
\end{Theorem}
The theorem states that the closed-loop \eqref{cascade-u}-\eqref{cascade-v-initial2} is exponentially stable, so is the original system \eqref{original-u}-\eqref{initial-u}. 

Before proving the stability, we first discuss the well-posedness of PDE \eqref{cascade-u}-\eqref{cascade-v-initial2} with control	\eqref{controller-U1} and \eqref{controller-U2}. The target system \eqref{target-w}-\eqref{target-z-initial2} is a standard 2-D heat equation with nonhomogeneous boundary conditions  at $y=l$ and $y=-l$ and the boundary conditions are given by $z_i$, $i=1,2$ which have been solved explicitly, so \eqref{target-w}-\eqref{target-z-initial2} has a unique solution \cite{myint2007linear}. 
Since the transformations \eqref{trans-w}, \eqref{trans-z1} and \eqref{trans-z2}  are invertible  and bounded, the original system  \eqref{cascade-u}-\eqref{cascade-v-initial2} with control \eqref{controller-U1} and \eqref{controller-U2} has a unique solution as well.

The proof of the stability includes two steps. First, we prove the stability of the target system \eqref{target-w}-\eqref{target-z-initial2} by introducing a change of variable. Second, we show the norm equivalence between cascade system $(u, v_1,v_2)$ and the target system $(w, z_1, z_2)$. Combining the above two steps, we finally reach Theorem 1.

Introduce the change of variable as
\begin{align}\label{trans-m-w}
m(x,y,t)=w(x,y,t)-\frac{y+l}{2l}z_1(x,0,t)+\frac{y-l}{2l}z_2(x,0,t),
\end{align}
which transforms the target system \eqref{target-w}-\eqref{target-z-initial2} into
\begin{align}
&m_t=m_{xx}+m_{yy}+\frac{(y+l)}{2l}\partial_{xx}z_1(x,0,t)\nonumber\\&~~~~~~~
-\frac{(y-l)}{2l}\partial_{xx}z_2(x,0,t)-\frac{(y+l)}{2l}\partial_tz_1(x,0,t)\nonumber\\
&~~~~~~~+\frac{(y-l)}{2l}\partial_tz_2(x,0,t),\label{target-m}\\
&m(-L,y,t)=m(L,y,t)=0,\label{target-m-boun1}\\
&m(x,-l,t)=m(x,l,t)=0\label{target-m-boun2},\\
&m(x,y,0)=m_0(x,y),\label{target-m-initial}\\
&\partial_t z_1(x,s,t)=\partial_s z_1(x,s,t)\label{re-target-z1-pde},\\
&\partial_t z_2(x,s,t)=\partial_sz_2(x,s,t)\label{re-target-z2-pde},\\
&z_1(x,\tau,t)=0\label{re-target-z-boundary1},\\
&z_2(x,\tau,t)=0\label{re-target-z-boundary2},\\
&z_1(x,s,0)=z_{10}(x,s),\label{re-target-z-initial1}\\
&z_2(x,s,0)=z_{20}(x,s)\label{re-target-z-re-initial2}.
\end{align}
\begin{Proposition}\label{stability-target}
\rm	For any initial conditions $(w_0,z_{10},z_{20})\in L^2(\mathcal{D}_1)\times H^1(\mathcal{D}_2)\times H^1(\mathcal{D}_2)$, and the compatible conditions satisfying
\begin{align}
	&w_0(L,y)=0,~~w_0(-L,y)=0,\\
	&w_0(x,l)=z_{10}(x,0),~~w_0(x,-l)=z_{20}(x,0),\\
	&z_{10}(x,\tau)=0,~~z_{20}(x,\tau)=0,
\end{align}
the system  \eqref{target-w}-\eqref{target-z-initial2} admits a unique solution and the equilibrium $(w, z_1, z_2) \equiv 0$ is exponentially stable, in particular, there  exist constants $\alpha_2, \beta_2>0$, such that
\begin{align}
	&\Vert w(t)\Vert^2_{L^2}+\Vert z_1(t)\Vert^2_{H^1}+\Vert z_2(t)\Vert^2_{H^1}\nonumber\\
	\le&\alpha_2e^{-\beta_2t}(\Vert w_0\Vert^2_{L^2}+\Vert z_{10}\Vert^2_{H^1}+\Vert z_{20}\Vert^2_{H^1}).
\end{align}
\end{Proposition}
\begin{Proof} \rm
First, define a Lyapunov function as
\begin{align*}
	&V_1(t)=\int_{-L}^{L}\int_{-l}^{l}m^2dydx+\int_{-L}^{L}\int_{0}^{\tau}e^{bs}(z^2_1+(\partial_xz_1)^2+\nonumber\\
	&(\partial_sz_1)^2)dsdx+\int_{-L}^{L}\int_{0}^{\tau}e^{bs}(z^2_2+(\partial_xz_2)^2+(\partial_sz_2)^2)dsdx,\label{Lya-V1}
\end{align*}
with the constant $b>0$.

Then, applying the integration by parts,  combining \eqref{target-m}-\eqref{re-target-z-re-initial2} and Cauchy-Schwarz inequality, we get
\begin{align*}
	&\frac{d}{dt}\int_{-L}^{L}\int_{-l}^{l}m^2dydx
	\le-(2-\mu_1-\mu_2)\int_{-L}^{L}\int_{-l}^{l}m^2_xdydx\nonumber\\
	&-2\int_{-L}^{L}\int_{-l}^{l}m^2_ydydx+\frac{8l}{3\mu_1}\int_{-L}^{L}(\partial_xz_1(x,0,t))^2dx\nonumber\\
	&+\frac{8l}{3\mu_2}\int_{-L}^{L}(\partial_xz_2(x,0,t))^2dx+\frac{8l}{3\mu_3}\int_{-L}^{L}(\partial_tz_1(x,0,t))^2dx\nonumber\\
	&+\frac{8l}{3\mu_4}\int_{-L}^{L}(\partial_tz_2(x,0,t))^2dx+(\mu_3+\mu_4)\int_{-L}^{L}\int_{-l}^{l}m^2dydx,
\end{align*}
\begin{align*}
	&\frac{d}{dt}\int_{-L}^{L}\int_{0}^{\tau}e^{bs}(z^2_1+(\partial_xz_1)^2+(\partial_sz_1)^2)dsdx\nonumber\\
	&+\frac{d}{dt}\int_{-L}^{L}\int_{0}^{\tau}e^{bs}(z^2_2+(\partial_xz_2)^2+(\partial_sz_2)^2)dsdx\nonumber\\
	&=-\int_{-L}^{L}(z^2_1(x,0,t)+(\partial_xz_1(x,0,t))^2+(\partial_{t}z_1(x,0,t))^2)dx\nonumber\\
	&-b\int_{-L}^{L}\int_{0}^{\tau}e^{bs}(z^2_1+(\partial_xz_1)^2+(\partial_tz_1)^2)dsdx\nonumber\\
	&-\int_{-L}^{L}(z^2_2(x,0,t)+(\partial_xz_2(x,0,t))^2+(\partial_tz_2(x,0,t))^2)dx\nonumber\\
	&-b\int_{-L}^{L}\int_{0}^{\tau}e^{bs}(z^2_2+(\partial_xz_2)^2+(\partial_tz_2)^2)dsdx.
\end{align*}
Applying the Poincar\'{e} inequality, i.e.,
\begin{align}
	\frac{1}{4L^2}\int_{-L}^{L}\int_{-l}^{l}m^2dydx&\le\int_{-L}^{L}\int_{-l}^{l}m^2_xdydx,\\
	\frac{1}{4l^2}\int_{-L}^{L}\int_{-l}^{l}m^2dydx&\le\int_{-L}^{L}\int_{-l}^{l}m^2_ydydx,
\end{align} 
we infer that
\begin{align}
	&\dot{V}_1(t)\le
	-b\int_{-L}^{L}\int_{0}^{\tau}e^{bs}(z^2_1+(\partial_tz_1)^2+(\partial_xz_1)^2)dsdx\nonumber\\
	&-(\frac{2-\mu_1-\mu_2}{4L^2}+\frac{1}{2l^2}-\mu_3-\mu_4)\int_{-L}^{L}\int_{-l}^{l}m^2dydx\nonumber\\
	&-b\int_{-L}^{L}\int_{0}^{\tau}e^{bs}(z^2_2+(\partial_tz_2)^2+(\partial_xz_2)^2)dsdx\nonumber\\
	&-\left(1-\frac{8l}{3\mu_1}\right)\int_{-L}^{L}(\partial_xz_1(x,0,t))^2dx\nonumber\\
	&-\left(1-\frac{8l}{3\mu_2}\right)\int_{-L}^{L}(\partial_xz_2(x,0,t))^2dx\nonumber\\
	&-\left(1-\frac{8l}{3\mu_3}\right)\int_{-L}^{L}(\partial_tz_1(x,0,t))^2dx\nonumber\\
	&-\left(1-\frac{8l}{3\mu_4}\right)\int_{-L}^{L}(\partial_tz_2(x,0,t))^2dx\nonumber\\
	&-\int_{-L}^{L}z^2_1(x,0,t)dx-\int_{-L}^{L}z^2_2(x,0,t)dx,
\end{align}
where $\frac{2-\mu_1-\mu_2}{4L^2}+\frac{1}{2l^2}-\mu_3-\mu_4>0
$, $1-\frac{8l}{3\mu_i}>0$, $i=1,2,3,4$.
Then, there exist constants $\alpha_3>0$ and $\beta_3>0$, such that
\begin{align}
	V_1\le \alpha_3e^{-\beta_3t} V_1(0).
\end{align}
To prove the stability of target system \eqref{target-w}-\eqref{target-z-initial2}, we define 
\begin{align*}
	&V_2(t)=\int_{-L}^{L}\int_{-l}^{l}w^2dydx+\int_{-L}^{L}\int_{0}^{\tau}(z^2_1+(\partial_xz_1)^2+\nonumber\\
	&(\partial_sz_1)^2)dsdx+\int_{-L}^{L}\int_{0}^{\tau}(z^2_2+(\partial_xz_2)^2+(\partial_sz_2)^2)dsdx.\label{Lya-V2}
\end{align*}
By a lengthy calculation, one can get that
\begin{align}
	\alpha_4 V_2\le V_1\le \beta_4 V_2,
\end{align}
where $\alpha_4=\frac{1}{4}$, $\beta_4=3+e^{b\tau}$, which gives the norm equivalence between $V_1$ and $V_2$. Hence,
\begin{align}
	V_2\le \alpha_2e^{-\beta_2t} V_2(0).
\end{align}
This proposition has been proved.
\end{Proof}

Before proceeding, we present the following definition of the Fourier series.
\begin{Definition}\label{non-fourier-series}\rm
According to \cite{haberman2012applied}, (P296),  it is known that   $\left\{\varphi_n(x)|n\in N^+\right\}$ and $\left\{\phi_m(y)|m\in N^+\right\}$ defined in \eqref{varphi} and \eqref{phi}  are orthogonal eigenfunctions, respectively.  For functions $h(y)\in L^2([-l,l])$,  $f(x,y)\in L^2(\mathcal{D}_1)$, and $g(x,s)\in L^2(\mathcal{D}_2)$, it holds \cite{doi:10.1080/00207179.2017.1286693}
\begin{align}
	h(y)&=\sum_{m=1}^{\infty}h_m\phi_m(y),\\
	f(x,y)&=\sum_{n=1}^{\infty}\sum_{m=1}^{\infty}a_{n,m}(f)\varphi_n(x)\phi_m(y),\\
	g(x,s)&=\sum_{n=1}^{\infty}g_n(s)\varphi_n(x),\label{fourier-g}
\end{align}
where the Fourier coefficients
\begin{align}
	h_m&=\frac{1}{l}\int_{-l}^{l}h(\xi)\phi_m(\xi)d\xi,\\
	a_{n,m}(f)&=\frac{1}{Ll}\int_{-L}^{L}\int_{-l}^{l}f(\theta,\xi)\varphi_n(\theta)\phi_m(\xi)d\xi d\theta,\\
	g_n(s)&=\frac{1}{L}\int_{-L}^{L}g(\theta,s)\varphi_n(\theta)d\theta.
\end{align}
From the Parseval's theorem, we know   $\Vert h\Vert_{L^2}^2=l\sum_{m=1}^{\infty}h^2_m$,  $\Vert f\Vert^2_{L^2}=Ll\sum_{n=1}^{\infty}\sum_{m=1}^{\infty}a^2_{n,m}(f)$, and $\Vert g(x,s)\Vert^2_{L^2}=L\sum_{n=1}^{\infty}\int_{0}^{\tau}g^2_n(s)ds$. 

\end{Definition}

The following lemmas will be used in the proof of the norm equivalence between the original system and the target system.
\begin{Lemma}\label{inequ-p}\rm
For function $p(l,\xi)$, $\xi\in[-l,l]$, there exist positive constants $M_i$, $i=1,2,3,4$, such that $\Vert p(l,\xi)\Vert^2_{L^2}\le M_1$, $\Vert p_\xi(l,\xi)\Vert^2_{L^2}\le M_2$, $\Vert p_{\xi\xi}(l,\xi)\Vert^2_{L^2}\le M_3$, $p^2_{\xi}(l,\xi)|_{\xi=l}=M_4$.
\end{Lemma}
One can use Parseval's theorem and Fourier series to prove lemma \ref{inequ-p}. Here the proof is omitted and we refer the reader to  \cite{krstic2009control}.
\begin{Lemma}\label{boundness-transformation}\rm
For  variables $(x,\theta,s,\xi)\in\Gamma_2$,  and  functions $f\in L^2(\mathcal{D}_1)$  and $g\in L^2(\mathcal{D}_2)$, there exist positive constants $A_i$, $B_{i}$, $i=1,2$, such that
\begin{align}
	&\int_{-L}^{L}\int_{0}^{\tau}\left|\int_{-L}^{L}\int_{-l}^{l}\gamma_i(x,\theta,s,\xi)f(\theta,\xi)d\xi d\theta\right|^2dsdx\nonumber\\
	&\le A_{i}\Vert f\Vert^2_{L^2}\label{boundness--gamma-trans1},\\
	&\int_{-L}^{L}\int_{0}^{\tau}\left|\int_{-L}^{L}\int_{0}^{s}\partial_\xi\gamma_i(x,\theta,s-r,\pm l)g(\theta,r)dr d\theta\right|^2dsdx\nonumber\\
	&\le B_i\Vert g\Vert^2_{L^2}\label{boundness-M}.
\end{align}
\end{Lemma}
\begin{Proof} \rm
According to Definition \ref{non-fourier-series}, we know that $p_{1m}$ defined in \eqref{p-fourier} is the Fourier coefficient of $p(l,\xi)$, i.e.,
\begin{align}\label{fourier-p}
	p(l,\xi)=\sum_{m=1}^{\infty}p_{1m}\phi_m(\xi).
\end{align}
As $i=1$, applying the Cauchy-Schwarz inequality, then combining \eqref{solution-gamma1}, and Lemma \ref{inequ-p},
it gets
\begin{align*}
	&\int_{-L}^{L}\int_{0}^{\tau}\left|\int_{-L}^{L}\int_{-l}^{l}\gamma_1(x,\theta,s,\xi)f(\theta,\xi)d\xi d\theta\right|^2dsdx\\
	\le&e^{2\lambda \tau}l^2\int_{-L}^{L}\int_{0}^{\tau}\sum_{m=1}^{\infty}e^{\frac{-m^2\pi^2s}{2l^2}}p^2_{1m}\sum_{m=1}^{\infty}\left(\sum_{n=1}^{\infty}e^{\frac{-n^2\pi^2s}{4L^2}}\right.\\
	&\left.\cdot\varphi_n(x)a_{n,m}(f)\right)^2dsdx\\
	\le&e^{2\lambda \tau}l^2L\int_{0}^{\tau}\sum_{m=1}^{\infty}p^2_{1m}\sum_{m=1}^{\infty}\sum_{n=1}^{\infty}a^2_{n,m}(f)ds\nonumber\\
	\le&\tau e^{2\lambda\tau}\Vert p(l,\xi)\Vert^2_{L^2}\Vert f\Vert^2_{L^2}=A_1\Vert f\Vert^2_{L^2}.
\end{align*}
Employing a similar approach, one can get the same result as $i=2$ for inequality \eqref{boundness--gamma-trans1}.

Before proceeding, based on \eqref{fourier-p}, the Parseval's theorem, and the definition of $\phi_m(\xi)$ which is shown in \eqref{phi}, we can get
\begin{align*}
	p_{\xi}(l,\xi)=\sum_{m=1}^{\infty}p_{1m}\phi'_m(\xi),\Vert p_{\xi}(l,\xi)\Vert^2_{L^2}=l\sum_{m=1}^{\infty}\frac{m^2\pi^2}{4l^2}p^2_{1m}.
\end{align*}
For the inequality \eqref{boundness-M}, since the proof process of four different cases of inequality \eqref{boundness-M} is similar, we just show the case as the kernel is $\partial_{\xi}\gamma_1(x,\theta,s-r,l)$.

Based on \eqref{solution-gamma1}, using the Cauchy-Schwarz inequality, Definition  \ref{non-fourier-series} and Lemma \ref{inequ-p}, one gets
\begin{align*}
	&\int_{-L}^{L}\int_{0}^{\tau}\left|\int_{-L}^{L}\int_{0}^{s}\partial_\xi\gamma_1(x,\theta,s-r,l)g(\theta,r)dr d\theta\right|^2dsdx\\
	\le& e^{2\lambda\tau}\int_{0}^{\tau}\int_{0}^{s}\sum_{m=1}^{\infty}e^{\frac{-m^2\pi^2(s-r)}{2l^2}}\sum_{m=1}^{\infty}\frac{m^2\pi^2}{4l^2}p^2_{1m}dr\\
	&\cdot\int_{-L}^{L}\int_{0}^{s}\left(\sum_{n=1}^{\infty}e^{\frac{-n^2\pi^2(s-r)}{2L^2}}\varphi_n(x)\frac{1}{L}\int_{-L}^{L}\varphi_n(\theta)\right.\\
	&\left.\cdot g(\theta,r)d\theta\right)^2 dr dxds\\
	\le& e^{2\lambda\tau}\frac{L}{l}\Vert p_\xi(l,\xi)\Vert^2_{L^2}\sum_{m=1}^{\infty}\frac{2l^2}{m^2\pi^2}\int_{0}^{\tau}\int_{0}^{s}\sum_{n=1}^{\infty}e^{\frac{-n^2\pi^2(s-r)}{L^2}}\\
	&\cdot g^2_n(s)drds\\
	\le& e^{2\lambda\tau}L M_2\sum_{m=1}^{\infty}\frac{2l}{m^2\pi^2}\int_{0}^{\tau}\int_{0}^{\tau}\sum_{n=1}^{\infty}g^2_n(r)drds\\
	\le& e^{2\lambda\tau}\tau  M_2\sum_{m=1}^{\infty}\frac{2l}{m^2\pi^2}\Vert g\Vert^2_{L^2}\le B_{1}\Vert g\Vert^2_{L^2}.
\end{align*}
The inequality \eqref{boundness-M} gets proven, and the proof of this lemma is completed.
\end{Proof}

\begin{Remark}
\rm
Lemma \ref{boundness-transformation} implies that the backstepping transformations \eqref{trans-z1} and \eqref{trans-z2} are bounded, despite there being singularities in the kernel functions. 
\end{Remark}
To prove the norm equivalence between the input $v_i(x,s,t)$, $i=1,2$ and $z_i(x,s,t)$, $i=1,2$, we propose the following two lemmas.

\begin{Lemma}\label{gamma-s-ux-uy}\rm
For  variables $(x,\theta,s,\xi)\in\Gamma_2$  and  function $f(x,y)\in L^2(\mathcal{D}_1)$, $g(x,s)\in L^2(\mathcal{D}_2)$, there exist positive constants $C_{ij}$, $D_{ij}$, $i=1, 2$, $j=1, 2$, such that
\begin{align}
	&\int_{-L}^{L}\int_{0}^{\tau}\left(\sum_{n=1}^{\infty}\sum_{m=1}^{\infty}\int_{-L}^{L}\int_{-l}^{l}\frac{(n\pi)^2}{4L^2}\right.\nonumber\\
	&\left.\gamma_i(x,\theta,s,\xi)f(\theta,\xi)d\xi d\theta\right)^2dsdx
	\le C_{1i}\Vert f_x\Vert^2_{L^2},\label{inequ-n-gamma}\\
	&\int_{-L}^{L}\int_{0}^{\tau}\left(\sum_{n=1}^{\infty}\sum_{m=1}^{\infty}\int_{-L}^{L}\int_{-l}^{l}\frac{(m\pi)^2}{4l^2}\right.\nonumber\\
	&\left.\gamma_i(x,\theta,s,\xi)f(\theta,\xi)d\xi d\theta\right)^2dsdx\le C_{2i}\Vert f_y\Vert^2_{L^2},\label{inequ-m-gamma}\\
	&\int_{-L}^{L}\int_{0}^{\tau}\left|\sum_{n=1}^{\infty}\sum_{m=1}^{\infty}\int_{-L}^{L}\int_{0}^{s}\frac{(n\pi)^2}{4L^2}\right.\label{inequ-M1}\\
	&\cdot\left.\partial_\xi\gamma_i(x,\theta,s-r,\pm l)g(\theta,r)dr d\theta\right|^2dsdx\le D_{1i}\Vert g_x\Vert^2_{L^2},\nonumber\\
	&\int_{-L}^{L}\int_{0}^{\tau}\left|\sum_{n=1}^{\infty}\sum_{m=1}^{\infty}\int_{-L}^{L}\int_{0}^{s}\frac{(m\pi)^2}{4l^2}\right.\label{inequ-M2}\\
	&\cdot\left.\partial_\xi\gamma_i(x,\theta,s-r,\pm l)g(\theta,r)dr d\theta\right|^2dsdx\le D_{2i}\Vert g\Vert^2_{L^2}.\nonumber
\end{align}
\end{Lemma}
The proof of this lemma is shown in Appendix \ref{proof-lemma-gamma-s-ux-uy}.
\begin{Lemma}\label{boundness-transformation-xs}\rm
For any variables $(x,\theta,s,\xi)\in\Gamma_2$,  and any functions $f\in L^2(\mathcal{D}_1)$,  $g\in L^2(\mathcal{D}_2)$, there exist positive constants  $E_{i}$, $F_{i}$, $G_{i}$, $H_{i}$, $K_{i}$, $i=1,2$, such that
\begin{align}
	&\int_{-L}^{L}\int_{0}^{\tau}\left|\int_{-L}^{L}\int_{-l}^{l}\partial_{s}\gamma_i(x,\theta,s,\xi)f(\theta,\xi)d\xi d\theta\right|^2dsdx\nonumber\\
	&\le E_{i}(\Vert f\Vert^2_{L^2}+\Vert f_x\Vert^2_{L^2}+\Vert f_y\Vert^2_{L^2})\label{boundness--gamma-partial-s},\\
	&\int_{-L}^{L}\int_{0}^{\tau}\left|\int_{-L}^{L}\int_{0}^{s}\partial_{\xi s}\gamma_i(x,\theta,s-r,\pm l)g(\theta,r)dr d\theta\right|^2dsdx\nonumber\\
	&\le F_{i}(\Vert g\Vert^2_{L^2}+\Vert g_x\Vert^2_{L^2})\label{boundness-M-s},\\
	&\int_{-L}^{L}\int_{0}^{\tau}\left|\int_{-L}^{L}\partial_{\xi }\gamma_i(x,\theta,0,\pm l)g(\theta,s) d\theta\right|^2dsdx\nonumber\\
	&\le G_{i}\Vert g\Vert^2_{L^2}\label{boundness-M-ss},\\
	&\int_{-L}^{L}\int_{0}^{\tau}\left|\int_{-L}^{L}\int_{-l}^{l}\partial_{x}\gamma_i(x,\theta,s,\xi)f(\theta,\xi)d\xi d\theta\right|^2dsdx\nonumber\\
	&\le H_{i}\Vert f_x\Vert^2_{L^2}\label{boundness--gamma-partial-x},\\
	&\int_{-L}^{L}\int_{0}^{\tau}\left|\int_{-L}^{L}\int_{0}^{s}\partial_{\xi x}\gamma_i(x,\theta,s-r,\pm l)g(\theta,r)dr d\theta\right|^2dsdx\nonumber\\
	&\le K_{i}\Vert g_x\Vert^2_{L^2}\label{boundness-M-x}.
\end{align}
\end{Lemma}
\begin{Proof} \rm	
For all inequalities, we just prove the case of $i=1$ ($\xi=l$ if it is in the inequality), the other cases can be proved in the same way. 

First, for the inequality \eqref{boundness--gamma-partial-s}, using \eqref{solution-gamma1}, the fundamental inequality, \eqref{boundness--gamma-trans1} in Lemma \ref{boundness-transformation} and \eqref{inequ-n-gamma}, \eqref{inequ-m-gamma} in Lemma \ref{gamma-s-ux-uy}, it gets
\begin{align*}
	&\int_{-L}^{L}\int_{0}^{\tau}\left|\int_{-L}^{L}\int_{-l}^{l}\partial_s\gamma_1(x,\theta,s,\xi)f(\theta,\xi)d\xi d\theta\right|^2dsdx\\
	\le&3\lambda^2\int_{-L}^{L}\int_{0}^{\tau}\left|\int_{-L}^{L}\int_{-l}^{l}\gamma_1(x,\theta,s,\xi)f(\theta,\xi)d\xi d\theta\right|^2dsdx\\
	&+3\int_{-L}^{L}\int_{0}^{\tau}\left|\sum_{n=1}^{\infty}\sum_{m=1}^{\infty}\int_{-L}^{L}\int_{-l}^{l}\frac{(n\pi)^2}{4L^2}\gamma_1(x,\theta,s,\xi)\right.\\
	&~~~~~~~~~~~~~~~~~~~~~~\cdot\left.f(\theta,\xi)d\xi d\theta\right|^2dsdx\\
	&+3\int_{-L}^{L}\int_{0}^{\tau}\left|\sum_{n=1}^{\infty}\sum_{m=1}^{\infty}\int_{-L}^{L}\int_{-l}^{l}\frac{(m\pi)^2}{4l^2}\gamma_1(x,\theta,s,\xi)\right.\\
	&~~~~~~~~~~~~~~~~~~~~~~\cdot\left.f(\theta,\xi)d\xi d\theta\right|^2dsdx\\
	&\le3(\lambda^2 A_1\Vert f\Vert^2_{L^2}+C_{11}\Vert f_x\Vert^2_{L^2}+C_{21}\Vert f_y\Vert^2_{L^2}).
\end{align*}

For the next inequality, employing \eqref{solution-gamma1}, the fundamental inequality, \eqref{boundness-M} in Lemma \ref{boundness-transformation} and \eqref{inequ-M1}, \eqref{inequ-M2} in Lemma \ref{gamma-s-ux-uy}, one gets
\begin{align*}
	&\int_{-L}^{L}\int_{0}^{\tau}\left|\int_{-L}^{L}\int_{0}^{s}\partial_{\xi s}\gamma_1(x,\theta,s-r, l)g(\theta,r)dr d\theta\right|^2dsdx\\
	&\le3\lambda^2\int_{-L}^{L}\int_{0}^{\tau}\left|\int_{-L}^{L}\int_{0}^{s}\partial_{\xi }\gamma_1(x,\theta,s-r, l)g(\theta,r)dr d\theta\right|^2\\
	&~~~~~~~~~~~~~~~~~~~~~~~ dsdx\\
	&+3\int_{-L}^{L}\int_{0}^{\tau}\left|\sum_{n=1}^{\infty}\sum_{m=1}^{\infty}\int_{-L}^{L}\int_{0}^{s}\frac{(m\pi)^2}{4l^2}\partial_{\xi }\gamma_1(x,\theta,s-r, l)\right.\\
	&~~~~~~~~~~~~~~~~~~~~~~\cdot\left.g(\theta,r)dr d\theta\right|^2dsdx\\
	&\le3(\lambda^2 B_{1}\Vert g\Vert^2_{L^2}+D_{11}\Vert g_x\Vert^2_{L^2}+D_{21}\Vert g\Vert^2_{L^2}).
\end{align*}

For inequality \eqref{boundness-M-ss}, applying \eqref{solution-gamma1},  the Cauchy-Schwarz inequality, \eqref{initial-gamma1} and Lemma \ref{inequ-p}, one can infer that
\begin{align*}
	&\int_{-L}^{L}\int_{0}^{\tau}\left|\int_{-L}^{L}\partial_{\xi }\gamma_1(x,\theta,0, l)g(\theta,s) d\theta\right|^2dsdx\\
	=&p^2_{\xi}(l,\xi)|_{\xi=l}\int_{-L}^{L}\int_{0}^{\tau}\left(\sum_{n=1}^{\infty}g_n(s)\varphi_n(x)\right)^2dsdx\\
	=&p^2_{\xi}(l,\xi)|_{\xi=l}\Vert g\Vert^2_{L^2}\le M_4\Vert g\Vert^2_{L^2}=G_1\Vert g\Vert^2_{L^2}.
\end{align*}

Combining inequality \eqref{inequ-n-gamma}, the equation $\int_{-L}^{L}(\varphi'(x))^2dx=\frac{n^2\pi^2}{4L^2}\int_{-L}^{L}\varphi^2_n(x)dx=\frac{n^2\pi^2}{4L^2}L$, and through the proof process similar to the proof of Lemma \ref{gamma-s-ux-uy}, one gets
\begin{align*}
	&\int_{-L}^{L}\int_{0}^{\tau}\left|\int_{-L}^{L}\int_{-l}^{l}\partial_{x}\gamma_1(x,\theta,s,\xi)f(\theta,\xi)d\xi d\theta\right|^2dsdx\nonumber\\
	\le&e^{2\lambda \tau}l^2\int_{-L}^{L}\int_{0}^{\tau}\sum_{m=1}^{\infty}e^{\frac{-m^2\pi^2s}{2l^2}}p^2_{1m}\sum_{m=1}^{\infty}\left(\sum_{n=1}^{\infty}e^{\frac{-n^2\pi^2s}{4L^2}}\right.\nonumber\\
	&\left.\cdot\varphi'_n(x)a_{n,m}(f)\right)^2dsdx\nonumber\\
	\le&e^{2\lambda \tau}l^2L\int_{0}^{\tau}\sum_{m=1}^{\infty}p^2_{1m}\sum_{m=1}^{\infty}\sum_{n=1}^{\infty}e^{\frac{-n^2\pi^2s}{2L^2}} \frac{n^2\pi^2}{4L^2}a^2_{n,m}(f)ds\nonumber\\
	\le&H_1\Vert f_x\Vert^2_{L^2}.
\end{align*}
For the inequality \eqref{boundness-M-x}, applying \eqref{norm_g_x} in Appendix \ref{proof-lemma-gamma-s-ux-uy} and the approach similar to the proof of the inequality \eqref{inequ-M1}, we can get
\begin{align*}
	&\int_{-L}^{L}\int_{0}^{\tau}\left|\int_{-L}^{L}\int_{0}^{s}\partial_{\xi x}\gamma_1(x,\theta,s-r,l)g(\theta,r)dr d\theta\right|^2dsdx\nonumber\\
	=&\int_{-L}^{L}\int_{0}^{\tau}\left(\int_{-L}^{L}\int_{0}^{s}\sum_{n=1}^{\infty}\sum_{m=1}^{\infty}e^{\left(-\frac{(n\pi)^2}{4L^2}-\frac{(m\pi)^2}{4l^2}+\lambda\right)(s-r)}\right.\nonumber\\
	&\cdot (-1)^m\frac{m\pi}{2l}p_{1m}\frac{1}{L}\varphi'_n(x)\varphi_n(\theta)\left. g(\theta,r)dr d\theta\right)^2dsdx\nonumber\\
	=&\int_{-L}^{L}\int_{0}^{\tau}\left(\int_{-L}^{L}\int_{0}^{s}\sum_{n=1}^{\infty}\sum_{m=1}^{\infty}e^{\left(-\frac{(n\pi)^2}{4L^2}-\frac{(m\pi)^2}{4l^2}+\lambda\right)(s-r)}\right.\nonumber\\
	&\cdot(-1)^m\frac{m\pi}{2l}p_{1m}\frac{1}{L}\frac{n\pi}{2L}\varphi_n(x)\varphi_n(\theta)\left. g(\theta,r)dr d\theta\right)^2dsdx\nonumber\\
	\le&L\int_{0}^{\tau}\sum_{m=1}^{\infty}\frac{m^2\pi^2}{4l^2}p^2_{1m}\sum_{m=1}^{\infty}\int_{0}^{s}e^{2(\frac{-m^2\pi^2}{4l^2}+\lambda)(s-r)}\nonumber\\
	&\cdot\left.\sum_{n=1}^{\infty}\int_{0}^{s}\frac{n^2\pi^2}{4L^2} g^2_n(r) dr\right)^2 ds\nonumber\\
	\le&\tau e^{2\lambda\tau} \frac{L}{l}\Vert p_{\xi}(l,\xi)\Vert^2_{L^2}\sum_{m=1}^{\infty}\frac{2l^2}{m^2\pi^2}\sum_{n=1}^{\infty}\int_{0}^{\tau}\frac{n^2\pi^2}{4L^2} g^2_n(r) dr\nonumber\\
	\le&K_1\Vert g_x\Vert^2_{L^2}.
\end{align*}
This inequality gets proven.

So far, we have proved this lemma.
\end{Proof}

The following proposition states that the original system \eqref{cascade-u}-\eqref{cascade-v-initial2} and the target system \eqref{target-w}-\eqref{target-z-initial2} are equivalent in the sense of norm.

\begin{Proposition}\rm
There exist positive constants $\alpha_5$ and $\beta_5$, such that
\begin{align}
	&\alpha_5(\Vert u\Vert^2_{L^2}+\Vert v_1\Vert^2_{H^1}+\Vert v_2\Vert^2_{H^1})\nonumber\\
	\le&\Vert w\Vert^2_{L^2}+\Vert z_1\Vert^2_{H^1}+\Vert z_2\Vert^2_{H^1}\nonumber \\
	\le&\beta_5(\Vert u\Vert^2_{L^2}+\Vert v_1\Vert^2_{H^1}+\Vert v_2\Vert^2_{H^1}).\label{original-target}
\end{align}
According to the transformations \eqref{trans-w}-\eqref{trans-z2}, utilizing fundamental inequality, Cauchy-Schwarz inequality, and Lemma \ref{inequ-p}, \ref{boundness-transformation}, \ref{boundness-transformation-xs}, we can prove this lemma. Here is omitted.
\end{Proposition}

\begin{figure*}[htbp]
\centering
\subfloat[]{\epsfig{figure=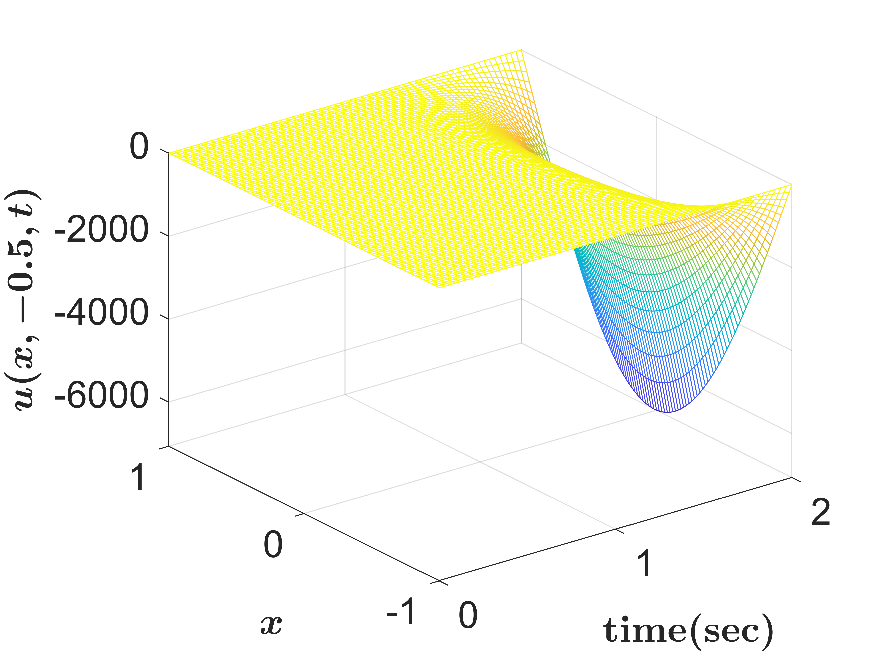,width=0.22\textwidth}}
\subfloat[]{\epsfig{figure=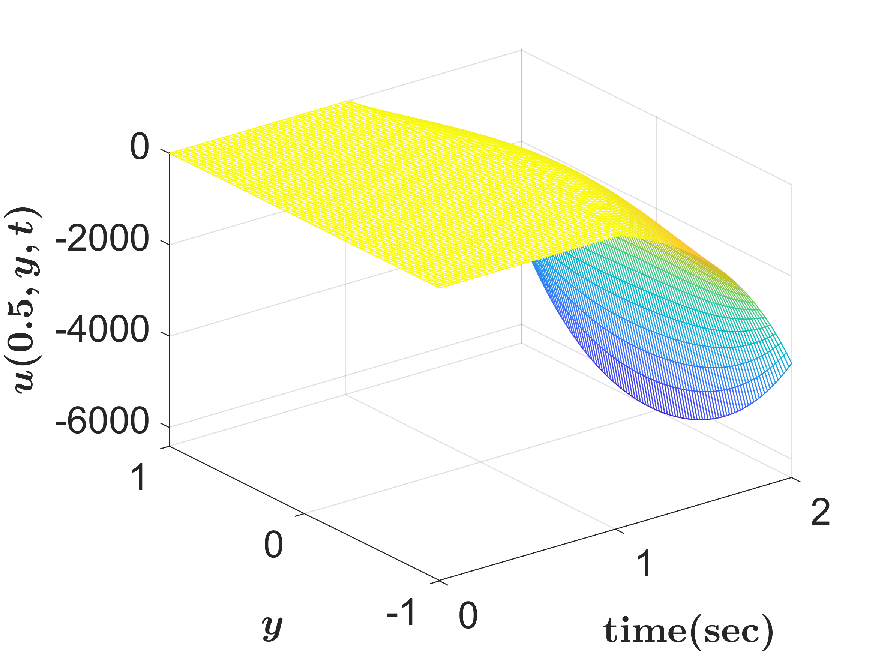,width=0.22\textwidth}}
\subfloat[]{\epsfig{figure=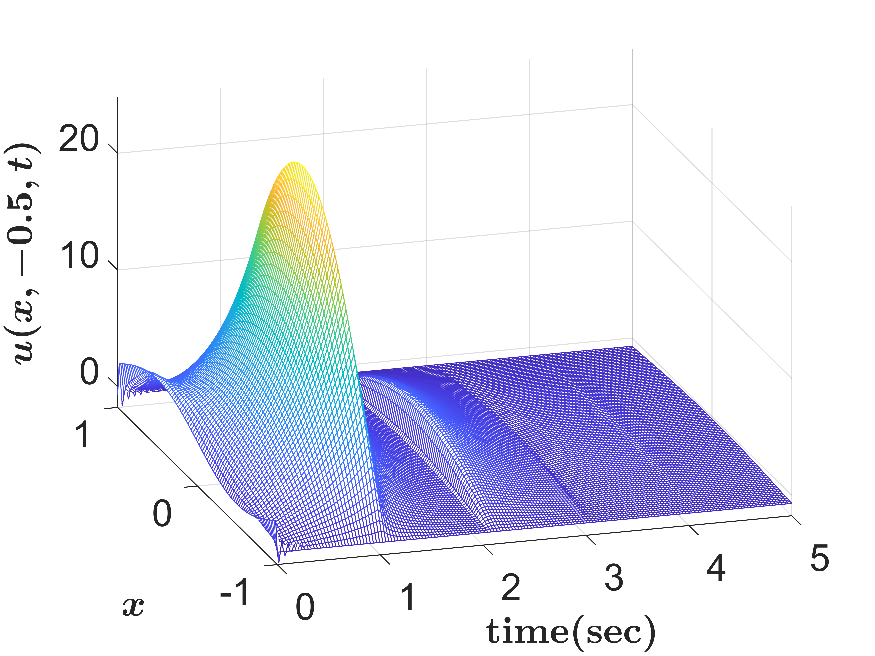,width=0.22\textwidth}}
\subfloat[]{\epsfig{figure=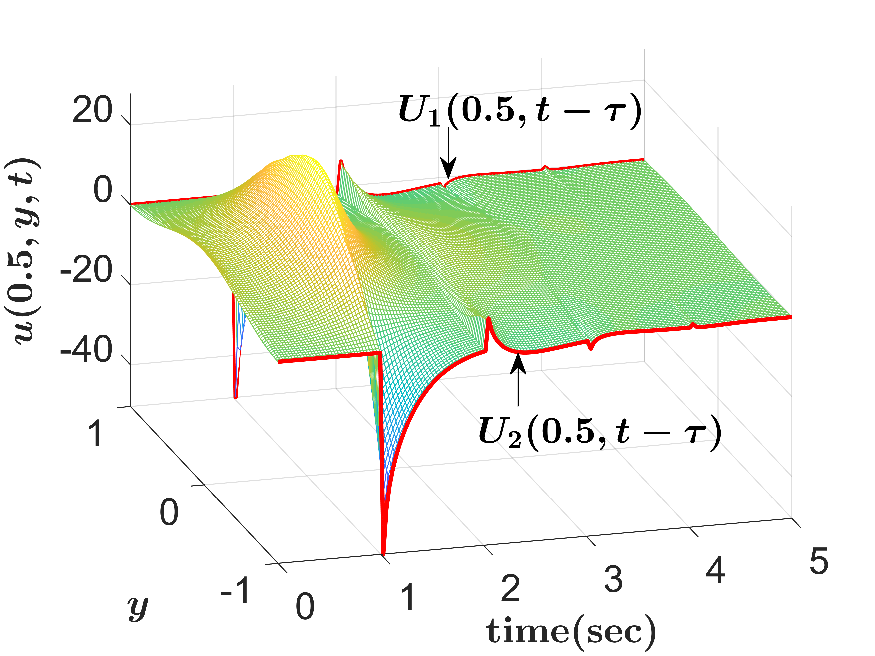,width=0.22\textwidth}}
\caption{The dynamics of the state with $\lambda=7$. (a)  $u(x,-0.5,t)$ without delay compensation. (b)  $u(0.5,y,t)$ without delay compensation. (c)  $u(x, -0.5, t)$ with delay compensation. (d) $u(0.5,y, t)$ with delay compensation. }\label{fig:u}
\end{figure*}
\begin{figure}[htbp]
\centering
\subfloat[]{\epsfig{figure=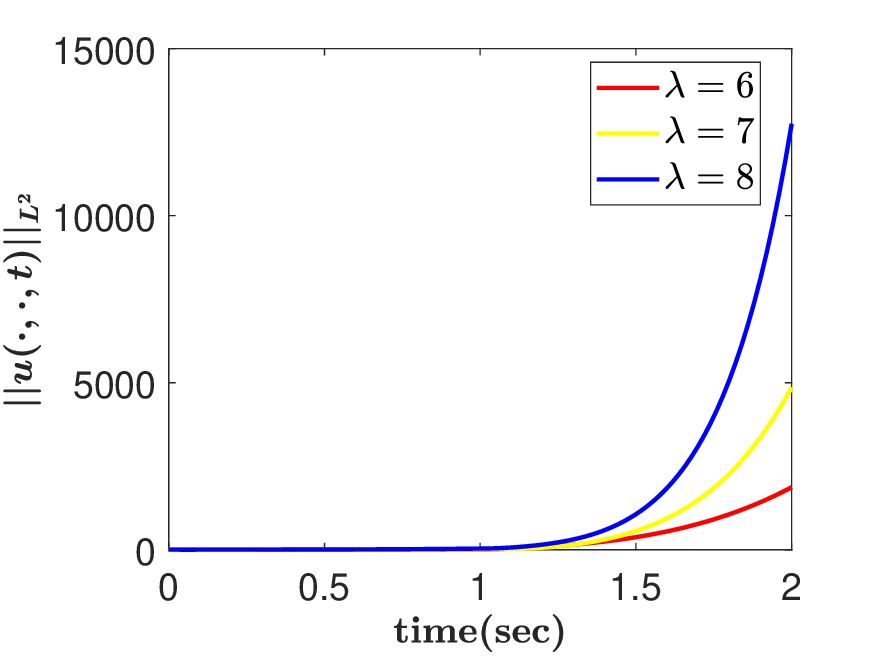,width=0.22\textwidth}}
\subfloat[]{\epsfig{figure=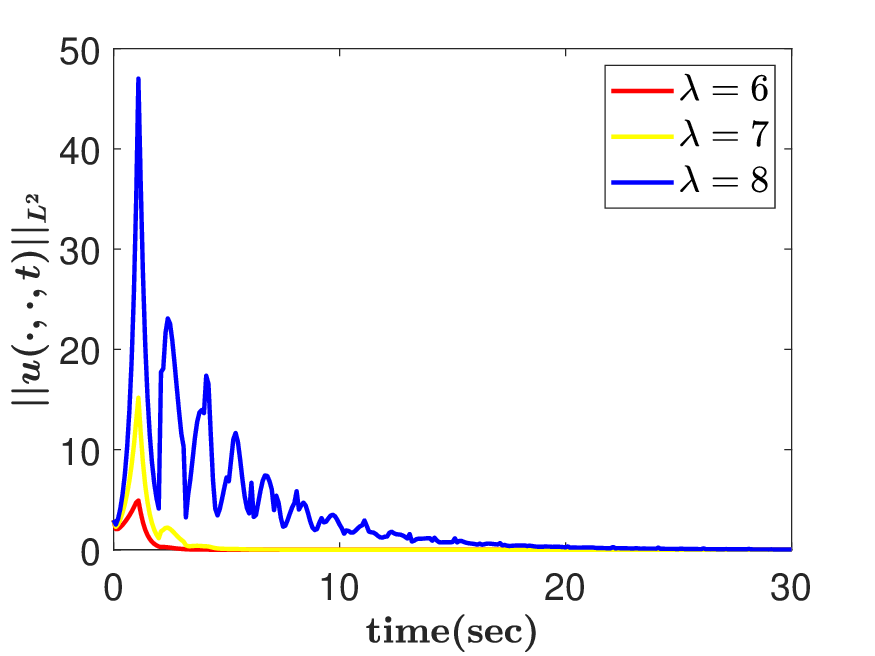,width=0.22\textwidth}}
\caption{The norm $||u(x,y,t)||_{L^2}$. (a) Without delay compensation. (b) With delay compensation.}\label{fig:L2}.
\end{figure}
\begin{figure}[htbp]
\centering
\subfloat[]{\epsfig{figure=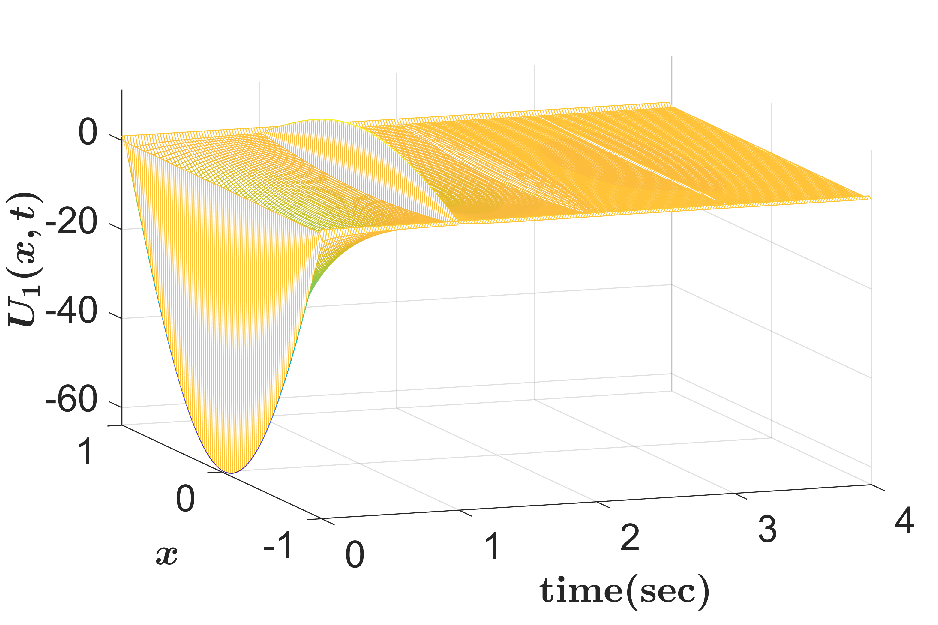,width=0.22\textwidth}}
\subfloat[]{\epsfig{figure=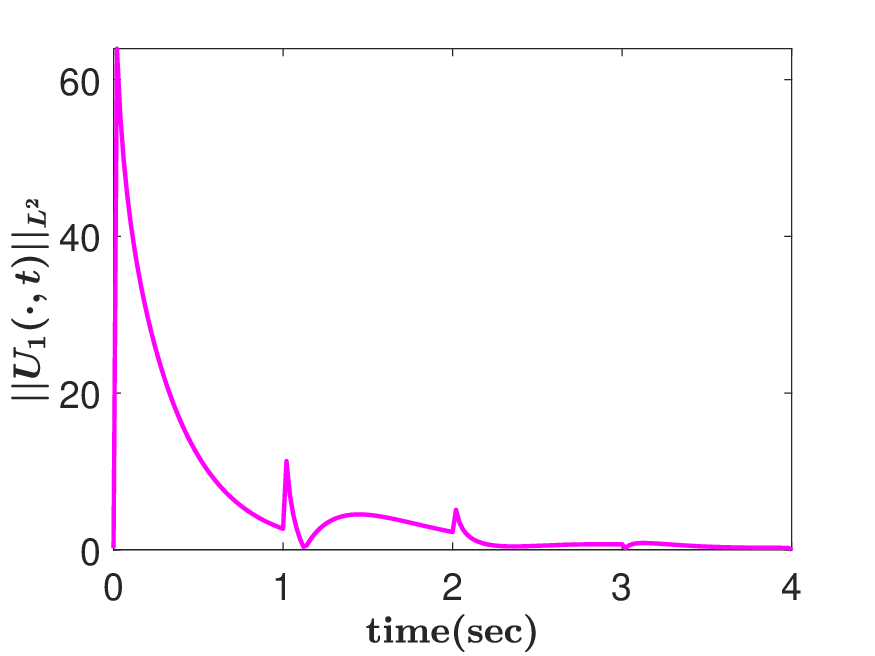,width=0.22\textwidth}}
\caption{The actuator $U_1(\cdot,t)$ with $\lambda=7$. (a) Dynamics. (b)  $\Vert U_1(\cdot,t)\Vert_{L^2}$.}\label{fig:U}
\end{figure}

\section{Simulation}
In this section, we provide an example to illustrate the effectiveness of the proposed delay-compensated bilateral control for the PDE system. In the numerical example, the domain is defined as $\mathcal{D}_1=[-1,1]\times[-1,1]$, and choose time delay $\tau=1s$, initial condition $u_0=2(\sin(\pi x)+1)(\cos(\pi y)+1)$. We apply the finite difference method by discretizing the domain by step sizes   $\Delta x=0.02$, $\Delta y=0.02$, $\Delta s=0.01$ and $\Delta t=0.01$.

The dynamics for both compensated and uncompensated cases with $\lambda=7$ are shown in Figure \ref{fig:u}, respectively. To observe the effectiveness of the delay-compensated control actuators, the variation of norm under three different parameters is demonstrated. The norm of state diverges in the case without delay compensated control in Figure \ref{fig:L2} (a). Figure \ref{fig:L2} (b) illustrates the convergence rate of the state with delay compensation. The dynamics and norm of $U_1(x,t)$ are shown in Figure \ref{fig:U} when $\lambda=7$. 

\section{Conclusion}
This paper presents bilateral actuators for reaction-diffusion with arbitrarily large input delay in a rectangular domain. To design the bilateral delay-compensated actuators, we propose a novel set of backstepping transformations, which results in the kernel functions being expressed by double trigonometric series with singularities. The main difficulty is to prove the norm equivalence in the sense of  $H^1$ norm between the target system and the cascade system based on the kernel functions containing singularities in the rectangular domain.  To overcome this difficulty, the Cauchy-Schwarz inequality, the 2-D Fourier series, and the Parseval's theorem are applied in the proving.  Another difficulty is to prove the exponential stability of the target which demands designing a change of variables for target system. Finally, we apply the Lyapunov function to prove the exponential stability of the closed-loop system. Future research will consider control design problems with non-constant input delays, such as delays varying with spatial parameters or varying with time.


\bibliographystyle{elsarticle-num}        
\bibliography{reference}           



\begin{appendices}
\section{The proof of Lemma \ref{gamma-s-ux-uy}}\label{proof-lemma-gamma-s-ux-uy}
\begin{Proof} \rm
	Combining \eqref{solution-gamma1} and Definition \ref{non-fourier-series}, it infers that
	\begin{align*}
		&\int_{-L}^{L}\int_{0}^{\tau}\left(\sum_{n=1}^{\infty}\sum_{m=1}^{\infty}\int_{-L}^{L}\int_{-l}^{l}\frac{(n\pi)^2}{4L^2}\gamma_1(x,\theta,s,\xi)\right.\nonumber\\
		&\left.~~~~~~~~~~~~~~~\cdot f(\theta,\xi)d\xi d\theta\right)^2dsdx\\
		\le&e^{2\lambda\tau}\int_{-L}^{L}\int_{0}^{\tau}\left|\sum_{n=1}^{\infty}\sum_{m=1}^{\infty}\int_{-L}^{L}\int_{-l}^{l}\frac{(n\pi)^2}{4L^2}e^{-\frac{(n\pi)^2s}{4L^2}}\right.\nonumber\\
		&\cdot e^{\frac{-(m\pi)^2s}{4l^2}}\frac{1}{L}p_{1m}\varphi_n(x)\left.\varphi_n(\theta)\phi_m(\xi)f(\theta,\xi)d\xi d\theta\right|^2dsdx\\
		\le&e^{2\lambda\tau}l^2\int_{-L}^{L}\int_{0}^{\tau}\sum_{m=1}^{\infty}e^{-\frac{(m\pi)^2s}{2l^2}}p^2_{1m}\sum_{m=1}^{\infty}\left(\sum_{n=1}^{\infty}\frac{n\pi}{2L}e^{-\frac{(n\pi)^2s}{4L^2}}\right.\nonumber\\
		&\cdot\left.\varphi_n(x)\frac{n\pi}{2L}a_{n,m}(f)\right)^2dsdx\\
		\le&e^{2\lambda\tau}M_1lL\sum_{m=1}^{\infty}\sum_{n=1}^{\infty}\frac{n^2\pi^2}{4L^2}\int_0^{\tau}e^{-\frac{(n\pi)^2s}{2L^2}}\frac{n^2\pi^2}{4L^2}a^2_{n,m}(f)ds\\
		\le&e^{2\lambda\tau}M_1l L\sum_{m=1}^{\infty}\sum_{n=1}^{\infty}\frac{n^2\pi^2}{4L^2}a^2_{n,m}(f)\le C_1\Vert f_x\Vert^2_{L^2}.
	\end{align*}
	By a similar process, we can prove the inequality \eqref{inequ-m-gamma}.
	
	For \eqref{inequ-M1}, according to the Fourier series of $g(x,s)$ and $p(l,\xi)$ defined in \eqref{fourier-g} and \eqref{fourier-p}, we can get 
	\begin{align}\label{fourier-v-theat}
		g_{x}(x,s)&=\sum_{n=1}^{\infty}g_n(s)\varphi'_n(x),\\
		p_{\xi\xi}(l,\xi)&
		=-\sum_{m=1}^{\infty}\frac{m^2\pi^2}{4l^2}p_{1m}\phi_m(x),
	\end{align} 
	then, combining the Parseval's theorem, it gets
	\begin{align}
		\Vert g_x(x,s)\Vert^2_{L^2}&=L\sum_{n=1}^{\infty}\frac{(n\pi)^2}{4L^2}\int_{0}^{\tau}g^2_n(s)ds,\label{norm_g_x}\\
		\Vert p_{\xi\xi}(l,\xi)\Vert^2_{L^2}&=l\sum_{m=1}^{\infty}\frac{m^4\pi^4}{16l^4}p^2_{1m}.\label{norm-p-xi2}
	\end{align}
	Finally, utilizing the Cauchy-Schwarz inequality, Lemma \ref{inequ-p}, \eqref{norm_g_x} and \eqref{norm-p-xi2}, we can infer that
	\begin{align*}
		&\int_{-L}^{L}\int_{0}^{\tau}\left|\sum_{n=1}^{\infty}\sum_{m=1}^{\infty}\int_{-L}^{L}\int_{0}^{s}\frac{n^2\pi^2}{4L^2}\right.\nonumber\\
		&~~~~~~~~~~~~\cdot\left.\partial_\xi\gamma_1(x,\theta,s-r, l)g(\theta,r)dr d\theta\right|^2dsdx\\
		\le&\int_{-L}^{L}\int_{0}^{\tau}\sum_{m=1}^{\infty}\frac{m^4\pi^4}{16l^4}p^2_{1m}\sum_{m=1}^{\infty}\frac{4l^2}{m^2\pi^2}\left(\sum_{n=1}^{\infty}\int_{0}^{s}\frac{n^2\pi^2}{4L^2}\varphi_n(x)\right.\\
		&\cdot\left. e^{(\frac{-m^2\pi^2}{4l^2}+\frac{-n^2\pi^2}{4L^2}+\lambda)(s-r)} g_n(r) dr\right)^2  dsdx\\
		\le&e^{2\lambda\tau}Ll\sum_{m=1}^{\infty}\frac{m^4\pi^4}{16l^4}p^2_{1m}\sum_{m=1}^{\infty}\frac{4l}{m^2\pi^2}\int_{0}^{\tau}\sum_{n=1}^{\infty}\int_{0}^{s}\left(\frac{n\pi}{2L}g_n(r)\right)^2dr\\
		&\cdot \int_{0}^{s}\frac{n^2\pi^2}{4L^2}e^{-\frac{n^2\pi^2}{2L^2}(s-r)}dr ds\\
		\le&e^{2\lambda\tau}L\tau\Vert p_{\xi\xi}(l,\xi)\Vert^2_{L^2}\sum_{n=1}^{\infty}\int_{0}^{\tau}\left(\frac{n\pi}{2L}g_n(r)\right)^2dr\\
		\le&e^{2\lambda\tau}\tau M_3\Vert g_x\Vert^2_{L^2}.
	\end{align*}
	This inequality gets proved. By a similar approach, we can prove \eqref{inequ-M2}.
\end{Proof}
\end{appendices}
\end{document}